\documentclass[english]{amsart}
\usepackage[T1]{fontenc}
\usepackage[latin9]{inputenc}
\usepackage{amssymb}

\makeatletter
\numberwithin{equation}{section} 
\numberwithin{figure}{section} 
  \@ifundefined{theoremstyle}{\usepackage{amsthm}}{}
  \theoremstyle{plain}
  \newtheorem{thm}{Theorem}[section]
  \theoremstyle{remark}
  \newtheorem*{acknowledgement*}{Acknowledgement}

\makeatother

\usepackage{babel}

\begin{document}
\title[An Upper Bound on the critical density]{An Upper Bound on the critical density for Activated Random Walks on Euclidean Lattices}

\begin{abstract}
We show the critical density for activated random walks on Euclidean
lattices is at most one.
\end{abstract}

\author{Eric Shellef}

\date{October 6, 2008}

\maketitle

\section{Introduction}

Given a graph, the activated random walks (ARW) model starts with
an initial configuration in which each vertex is occupied by a finite
number of sleeping or active particles. Beginning with this initial
configuration, each active particle performs an independent, rate
one, random walk, while sleeping particles stay put. If a sleeping
particle occupies the same vertex as an active particle, it becomes
active immediately. Finally, active particles fall asleep independently
at a rate $\lambda>0$. We examine this model on $\mathbb{Z}^{d}$,
where we suppose that initially, each vertex contains an i.i.d. Poisson
number of active particles with expected value $\mu$. The Poisson
distribution plays no special role and could be replace by other distributions.

One obvious question on the long term behavior of the system, is whether
or not we have fixation, which by translation invariance is equivalent
to whether the number of active particles that visit the origin is
finite almost surely. In Theorem \ref{thm:Consider-the-activated}
below, we show that for $\mu>1$, we almost surely do not have fixation.

For this, we rely on the technical framework developed in \cite{Rolla}.
In this paper, the existence of the process is proved, and it is shown
that the probability of finitely describable events can be approximated
by finite systems. Let $\mathbb{P}^{\mu}$ be the probability measure
on the model described above, and let $\mathbb{P}_{M}^{\mu}$ be the
measure on the model with all particles outside of $B_{M}$ removed,
where $B_{M}=\left\{ x\in\mathbb{Z}^{d}:\|x\|\le M\right\} $. If
$A$ is an event measurable with respect to what happens in some finite
subset up to some time $t<\infty$, \[
\mathbb{P}^{\mu}(A)=\lim_{M\to\infty}\mathbb{P}_{M}^{\mu}(A).\]

The second tool from \cite{Rolla}, is a graphical representation
for systems with finitely many particles that has the desirable properties
of monotonicity and commutativity of certain parameters. Here is a
loose description. Let there be a universal clock that will ring with
the appropriate rate, fix some label for each particle, and let there
be an i.i.d. sequence of labels, independent of the clock. Also, at
each site, let there be an i.i.d. sequence of envelopes, each one
containing some instruction to be performed. When the clock rings
for the first time, it will ring for the particle indicated by the
first label in the sequence and at that moment this particle will
perform some action. If the particle is sleeping, nothing happens.
If the particle is active, it will open the first envelope at that
site, burn the envelope and perform the action written inside. The
instruction may be to jump to a specific neighbor, or to try to sleep.
Thus there are two types of envelopes: jump envelopes and sleep envelopes.
If the particle tries to sleep but there are other particles on the
same site, the envelope is burned anyway. This representation has
the activated random walk process as a natural marginal, and is described
formally for one dimension in the above reference.

The commutativity property says the following. Suppose for a given
realization of the process, the system fixates, that is, all particles
are passive for all large enough times (of course starting with finitely
many particles this happens a.s.). Then by changing the label sequence
and the universal clock, the system will stabilize at exactly the
same state, except that some particles may be permuted. Furthermore,
the amount of envelopes burned at each site is also preserved. So
the final state of the system is determined by the initial conditions
(positions and types of the particles) and the by sequences of envelopes.
The second property, monotonicity, states the following. Suppose for
some realization of the envelopes and initial conditions there is
fixation. Take a new configuration by deleting some particles on the
original one, changing the type of some particles from active to sleeping
and inserting some sleep envelopes at some sites\textquoteright{}
envelope sequences. Then for this new configuration the system also
stabilizes and the final number of envelopes that are burned at each
site (not counting the ones inserted) does not increase.

This framework was used to show that there is at most one phase transition
in the model for Euclidean lattices and that in one dimension, there
exists a phase transition for some $0<\mu_{c}\le1$.

We utilize these properties to modify the finite approximations of
the process to {}``sleepier'' ones, and show that when the density
$\mu$ is higher than one, the number of visits to the origin still
goes to infinity almost surely, hence there is no fixation.

\section{Result}

Here we settle one of the open problems posed in the concluding remarks
of the \cite{Rolla}, and prove that for Euclidean lattices of all
dimensions, $\mu_{c}\le1$. This is shown by using the below theorem
in conjunction with Theorem 1.1 in \cite{Rolla}.

\begin{thm}
\label{thm:Consider-the-activated}Consider the activated random walk
model on $\mathbb{Z}^{d}$. Then for any $\mu>1$ the system does
not fixate. 
\end{thm}
\begin{proof}
For $r\in\mathbb{N}$, let $A_{r}$ be the event that the origin is
visited by an active particle at least $r$ times before fixation.
Fixing $d$, and $\mu>1$, we prove the above by showing that for
all $r\in\mathbb{N}$, \begin{equation}
\lim_{M\to\infty}\mathbb{P}_{M}^{\mu}\left(A_{r}\right)=1.\label{limM,r}\end{equation}
 Let $n=n(M,\omega)$ be the number of particles in the system. Since
$\mu>1$, for some $\epsilon=\epsilon(\mu)$, \begin{equation}
\lim_{M\to\infty}\mathbb{P}\left[n>(1+\epsilon)\left|B_{M}\right|\right]=1.\label{manyParticles}\end{equation}
 Let $F=F(M)$ be the event in \eqref{manyParticles} that happens
with high probability.

For any finite $M$, we can use the monotonicity and commutativity
proven in \cite{Rolla} to reduce the model to the following IDLA-like
process.

First notice that regardless of $\lambda$ (the rate by which particles
fall asleep), if two particles or more occupy the same site, they
will all be active almost surely. Considering one particle, we are
thus assured it will continue its random walk at least until it has
reach an unoccupied site. Conditioning on $F$, we fix some order
on $N=\lceil\left(1+\epsilon\right)\left|B_{M}\right|\rceil$ randomly
chosen particles, and by adding sleep envelopes, we make the remaining
$n-N$ particles static so they won't interfere. We modify the label
sequence so that each of the $N$ particles in turn begins a random
walk from its initial location until reaching a site unoccupied by
other particles (possibly its starting location - in which case it
wouldn't move at all). Once reaching such a site, we force the particle
to remain there forever by inserting a sleep envelope for each of
its movement attempts . Let $\mathcal{P}$ denote the probability
measure on this {}``embedded'' Markov process. Let $V$ be the number
of particles that visit zero before stopping.

By the monotonicity, to prove \eqref{limM,r} it suffices to show
that $V$ grows linearly with $M$ with high probability. Formally,
\begin{equation}
\mathcal{P}\left[V>\frac{\epsilon}{4}M\ \Big|\ F(M)\right]\to1\label{visitsTo0}\end{equation}

To prove \eqref{visitsTo0} we use an idea from the original IDLA
paper \cite{IDLA}, which is even simpler to apply in our setting.

Let $\left\{ X_{i}\right\} _{i=1,\ldots,N}$ be the starting locations
of the $N$ particles. Since the density is i.i.d., for any $x\in B_{M},\mathbb{P}_{M}^{\mu}\left[X_{i}=x\right]=\left|B_{M}\right|^{-1}$.
Unlike the real model, we let the walks continue forever, but mark
(e.g. by coloring) the locations where they first visit an unmarked
vertex. We start with all vertices unmarked. Let each walk run in
turn (without stopping) and mark the first unmarked vertex it visits.
This may be the initial placement of the particle (and will be in
most cases). We call the component of marked vertices at each step
the cluster.

Let $W$ be the number of walks that visit $0$ before exiting $B_{M}$.

Let $L$ be the number of walks that visit $0$ before exiting $B_{M}$,
but after leaving the cluster (i.e. after stopping in the original
model).

Note that $W-L$ is the number of visits to zero of particles that
haven't left the cluster, and is at most equal to $V$. Thus we have 

\begin{eqnarray*}
\mathcal{P}[V<\frac{\epsilon}{4}M] & < & \mathcal{P}[W-L<\frac{\epsilon}{4}M]\\
 & \le & \mathcal{P}[W-\frac{\epsilon}{4}M\le a]+\mathcal{P}[L\ge a]\end{eqnarray*}
for any real $a$. We choose $a=(1-\frac{\epsilon}{2})E\left[W\right]$.
We bound the above terms by calculating the expected value of $M$
and $L$. Let $\tau_{0}$ be the first hitting time of $0$ of a random
walk, and let $\tau_{M}$ be the first exit time from $B_{M}$. 

\[
E[W]=\sum_{i=1}^{N}\sum_{x\in B_{M}}\mathcal{P}_{x}[\tau_{0}<\tau_{M}]P(X_{i}=x)=\frac{N}{\left|B_{M}\right|}\sum_{x\in B_{M}}\mathcal{P}_{x}[\tau_{0}<\tau_{M}]\]

$E\left[L\right]$ is hard to calculate, but note that each walk that
contributes to $L$ can be tied to the unique point at which it exits
the cluster. Thus, by the Markov property, if we start a random walk
from each vertex in $B_{M}$ and let $\hat{L}$ be be the number of
such walks that hit $0$ before exiting $B_{M}$, we have $\mathcal{P}[L\ge a]\le\mathcal{P}[\hat{L}\ge a]$.

Thus, $E[\hat{L}]\ge E[L]$, and we have\[
E[\hat{L}]=\sum_{x\in B_{M}}\mathcal{P}_{x}[\tau_{0}<\tau_{M}]\]

So we have $E[W]\ge(1+\epsilon)E[\hat{L}]$ which gives us for $\epsilon<\frac{1}{2}$
\[
\mathcal{P}[L\ge a]\le\mathcal{P}[\hat{L}\ge(1-\frac{\epsilon}{2})(1+\epsilon)E[\hat{L}]]\le\mathcal{P}[\hat{L}\ge(1+\frac{\epsilon}{4})E[\hat{L}]]\]

We can lowerbound $E[\hat{L}]$ by using the Green function identity:
\[
\mathcal{P}_{x}[\tau_{0}<\tau_{M}]=\frac{G_{M}(x,0)}{G_{M}(0,0)}\]
 where $G_{M}(a,b)$ is the average number of visits of a random starting
at $a$ to $b$ before leaving $B_{M}$. By symmetry of the Green
function we can write\[
E[\hat{L}]=\sum_{x\in B_{M}}\mathcal{P}_{x}[\tau_{0}<\tau_{M}]=G_{M}(0,0)^{-1}\sum_{x\in B_{M}}G_{M}(0,x)=G_{M}(0,0)^{-1}E_{0}[\tau_{M}].\]
By the optional stopping theorem with the martingale $\|X(t)\|^{2}-t$
we have $E_{0}[\tau_{M}]=M^{2}$. Second, $G_{M}(0,0)=M$ for the
line and is smaller for higher dimensions (e.g. by the monotonicity
law for electric networks). 

Since $E[W]\ge E[\hat{L}]\ge M$, we have $\mathcal{P}[W-\frac{\epsilon}{4}M\le a]\le\mathcal{P}[W\le(1-\frac{\epsilon}{4})E[W]]$.

Since $\hat{L}$ and $W$ are both sums of indicators, we can use
standard concentration inequalities, and the lower bound on $E[\hat{L}]$,
to show exponential decay in $N$ of $\mathcal{P}[V<\frac{\epsilon}{4}M]$
which proves \eqref{visitsTo0}, and we are done.
\end{proof}
\begin{acknowledgement*}
Thanks to Gady Kozma and Vladas Sidoravicius for telling me about
this problem.
\end{acknowledgement*}


\end{document}